\documentclass[12pt]{article}
\usepackage{amssymb,amsmath,indentfirst}
\title{
{\bf The calculation of the adjacency indices\\ of type $D$ singularities of a front}
\author{Vyacheslav D. Sedykh}}

\begin{document}

\date{}
\maketitle

\begin{abstract}
The formula for the calculation of the adjacency indices of types $D_\mu^\pm$ Legendrian monosingularities to multisingularities of types $A_{\mu_1}^{k_1}\dots A_{\mu_p}^{k_p}$ is obtained in the paper. Namely the number of strata of a given type $A_{\mu_1}^{k_1}\dots A_{\mu_p}^{k_p}$ in the discriminant of a real simple singularity of any of types $D_\mu^\pm$ is calculated.

\end{abstract}

We study the adjacencies of singularities of the front of a proper generic Legendrian map with simple stable singularities. For all the necessary definitions and facts from the theory of Legendrian singularities, see \cite{Arn96},\cite{SedARMJ},\cite{Sed-kniga-2021}.

\medskip
By Arnold's theorem on Legendrian singularities, simple stable germs of Legendrian maps are numbered by the symbols
$$
A_{\mu}\,(\mu=1,2,\dots),\quad D_{\mu}^\pm\,(\mu=4,5,\dots),\quad E_6,E_7,E_8
$$
up to the Legendrian equivalence. Namely, the fronts of such germs are diffeomorphic to discriminants of real simple function germs of the corresponding types. If the number $\mu$ is odd, then Legendrian singularities of types $D_{\mu}^{+}$ and $D_{\mu}^{-}$ are Legendrian equivalent. Their types are usually denoted by $D_{\mu}$.

Let us consider the unordered set of symbols that are the types of germs of a generic Legendrian map $f:L\rightarrow V$ at the preimages of a point $y\in V$. The formal commutative product $\mathcal{A}$ of these symbols is called the type of the multisingularity of the map $f$ at the point $y$ (or the type of a monosingularity if $y$ has only one preimage). If $f^{-1}(y)=\emptyset$, then $\mathcal{A}=\mathbf{1}$. The set $\mathcal{A}_f$ of points $y$ at which the map $f$ has a multisingularity of type $\mathcal{A}$ is a smooth submanifold in the space $V$.

Assume that $f$ has a multisingularity of type $\mathcal{B}$ at a point $y$ lying in the closure $\overline{\mathcal{A}_f}$ of the manifold $\mathcal{A}_f$. From the Looijenga's paper \cite{Loenga}, it follows that $y$ has a neighbourhood $O(y)\subset V$ such that its intersection $O(y)\cap\sigma$ with any connected component $\sigma$ of the manifold $\mathcal{A}_f$ such that $\overline{\sigma}\ni y$ is contractible. The number of the connected components of the intersection $O(y)\cap\mathcal{A}_f$ depends only on $\mathcal{A}$ and $\mathcal{B}$, is denoted by $J_{\mathcal A}(\mathcal{B})$ and is called the adjacency index of a multisingularity of type $\mathcal{B}$ to a multisingularity of type $\mathcal{A}$. The index $J_{\bf 1}(\mathcal{B})$ is equal to the number of the connected components of the complement to the front of the map $f$ in $O(y)$.

The adjacency indices of multisingularities are calculated via adjacency indices of mono\-singularities (see \cite[Proposition 2.5]{Sed2012}). The adjacency indices of monosingularities of types $A_\mu$ are calculated by the formula \cite[Theorem 2.6]{Sed2012}. The indices of all adjacencies of monosingularities of types $D_\mu^\pm$ and $E_\mu$ for $\mu\leq6$ have been listed in \cite[Theorems 2.8 and 2.9]{Sed2012}. The indices $J_{\bf 1}(D_\mu^\pm)$ and $J_{\bf 1}(E_\mu)$ for higher $\mu$ were calculated by Vassiliev \cite{Vas-discr}.

Legendrian monosingularities of types $D_\mu^\pm$ can only be adjacent to multisingularities of types $D_\nu^\pm A_{\mu_1}^{k_1}\dots A_{\mu_p}^{k_p}$ and $A_{\mu_1}^{k_1}\dots A_{\mu_p}^{k_p}$. The adjacency indices $J_{D_\nu^\pm A_{\mu_1}^{k_1}\dots A_{\mu_p}^{k_p}}(D_\mu^\pm)$ were calculated in \cite{Sed2024}. The adjacency indices $J_{A_{\mu_1}^{k_1}\dots A_{\mu_p}^{k_p}}(D_\mu^\pm)$ are calculated in the present paper. 

\medskip
{\bf Example.} Singular points of the front of a generic Legendrian map with simple stable singularities divide it into smooth pieces (connected components of the manifold of monosingularities of type $A_1$). The number of such pieces in a small neighbourhood of a monosingularity of type $D_{\mu}^\delta$ is calculated by the following formulas:
$$
J_{A_1}(D_{2k}^\delta)=\frac{2k(k^2-1)}{3}+k(k+3)\left[\frac{3-\delta}{4}\right],\quad J_{A_1}(D_{2k+1})=\frac{k(k+1)(4k+5)}{6}+k,
$$
where $k\geq2$ and $[\cdot]$ is the integer part of a number.

\medskip
Below we get the explicit formula for adjacency indices of Legendrian monosingularities of types $D_\mu^\pm$ to multisingularities of all types $A_{\mu_1}^{k_1}\dots A_{\mu_p}^{k_p}$. 

Let
$$
f:\mathbb{R}^{n-1}\to \mathbb{R}^n,\quad
(t,q)\mapsto (x,q,u),\quad x=-\frac{\partial S(t,q)}{\partial t},\quad u=S(t,q)-t\frac{\partial S(t,q)}{\partial t},
$$
where
$$
t=(t_1,\dots,t_m),\quad x=(x_1,\dots,x_m),\quad q=(q_{m+1},\dots,q_{n-1}),
$$
and $S=S(t,q)$ is a smooth function of $t$ depending on the parameter $q$. Then the map $f$ is Legendrian. If $m=1$ and the function $S$ is given by the formula
$
S=t_1^{\mu+1}+q_{\mu-1}t_1^{\mu-1}+...+q_2t_1^2,
$
where $\mu\geq1$, then the germ of $f$ at zero has type $A_\mu$.

We will further assume that $m=2$, $\delta=\pm1$ and
$$
S=t_1^2t_2+\delta t_2^{\mu-1}+q_{\mu-1}t_2^{\mu-2}+\dots+q_3t_2^2,\quad \mu\geq4.
$$
Then $f$ has a singularity of type $D_\mu^\delta$ at zero ($D_\mu^+$ if $\delta=+1$ and $D_\mu^-$ if $\delta=-1$). 

\medskip
{\bf Lemma 1.} {\it The preimages $(t,q)$ of a given point $(x,q,u)$ under the map $f$ are determined by the points $t=(t_1,t_2)$ such that:

{\rm1)} $t_2$ is a multiple real root of the polynomial
\begin{equation}
\delta t_2^{\mu}+q_{\mu-1}t_2^{\mu-1}+\dots+q_3t_2^3+x_2t_2^2-ut_2-\frac{x_1^2}{4};
\label{P1}
\end{equation}

{\rm2)} $t_1=-\frac{x_1}{2t_2}$ if $t_2\neq0$ and $t_1^2=-x_2$ if $t_2=0$.}

\medskip
This statement follows directly from the definitions.

\medskip
{\bf Lemma 2.} (\cite[Lemma 5.1]{Sed2012}) {\it Let
$$
H_i=H_i(t,q)=t_2^{i}\,\frac{\partial^{i+1}S}{\partial t_2^{i+1}}-(-1)^{i+1}(i+1)!\,t_1^2,\quad i=1,\dots,\mu-1.
$$
Then the germ of the map $f$ at a point $(t,q)\neq 0$ has the following type:

{\rm1)} $A_\nu,1\leq \nu\leq \mu-1$ if $H_1=\dots=H_{\nu-1}=0, H_\nu\neq0$;

{\rm2)} $A_3$ if $t_1=t_2=0,q_3\neq0$;

{\rm3)} $D_{\nu}^\pm$, $4\leq \nu\leq \mu-1$ if $t_1=t_2=q_3=\dots=q_{\nu-1}=0$ and $\pm q_\nu>0$.}

\medskip
Lemmas 1 and 2 imply: 

\medskip
{\bf Lemma 3.} {\it Let us fix a point $y=(x,q,u)$ and assume that the polynomial {\rm(\ref{P1})} have zero root of multiplicity $\kappa\geq2$. Consider the preimages $z=(t,q)$ of $y$ under the map $f$. 

{\rm1)} If $\kappa=2$ and $x_2>0$, then there are no preimages $z$, where $t_2=0$.

{\rm2)} If $\kappa=2$ and $x_2<0$, then there are exactly two preimages $z$, where $t_2=0$. They have the coordinate $t_1=\pm\sqrt{-x_2}$. The germs of $f$ at these preimages have type $A_1$.

{\rm3)} If $\kappa\geq3$, then there is exactly one preimage $z$, where $t_2=0$. It has the coordinate $t_1=0$. The germ of $f$ at this preimage has type $A_3$ if $\kappa=3$ and $D_{\kappa}^{{\rm sign}\, q_{\kappa}}$ if $\kappa\geq4$.}

\medskip
By Lemma 1, it follows that for $x_1\neq0$ each polynomial of the form (\ref{P1}) having multiple real roots defines a point $(x,q,u)$ on the front of the map $f$. Since the set of everywhere-positive polynomials is connected (in the space of polynomials of a given degree), Lemma 3 implies:

\medskip
{\bf Lemma 4.} {\it {\rm1)} Points $(x,q,u)$ for values $x_1$ of the same sign belong to the same connected component of the manifold of multisingularities of the map $f$ if and only if 
\begin{equation}
\begin{tabular}{c}
\mbox{the locations in the increasing order of real roots of the corresponding}\\ 
\mbox{polynomials {\rm(\ref{P1})} taking into account their multiplicities and signs are the same.}
\end{tabular} 
\label{cond1}
\end{equation}

{\rm2)} Let the intersection of a connected component of the manifold of multisingularities with the space $x_1\neq0$ is not empty. Then the intersection of this component with the hyperplane $x_1=0$ is not empty if and only if 
\begin{equation}
\begin{tabular}{c}
\mbox{the polynomial {\rm(\ref{P1})} corresponding to any point $(x,q,u)$ of the component}\\ 
\mbox{with the coordinate $x_1\neq0$ has a simple real root such that there are}\\ 
\mbox{no other real roots of the polynomial between them and zero.}
\end{tabular} 
\label{cond2}
\end{equation}
Moreover, for any point $(x,q,u)$ of this intersection either $u\neq0$ or $u=0,x_2>0$. 

{\rm3)} Points $(x,q,u)$ for values $x_1$ of opposite signs belong to the same connected component of the manifold of multisingularities if and only if both conditions {\rm(\ref{cond1})} and {\rm(\ref{cond2})} are satisfied.}

\medskip
Connected components of manifolds of multisingularities of the map $f$ are semialgebraic. Therefore from Lemma 4, we get the following 

\medskip
{\bf Theorem 1.} {\it Let $\mathcal{A}=A_{\mu_1}^{k_1}\dots A_{\mu_p}^{k_p}$. Then the connected components of the manifold $\mathcal{A}_f$ have the following types:

{\rm1)} the components lying in the space
$$
\Pi_1: x_1^2+u^2\neq0 \mbox{ or } x_1=u=0,x_2>0;
$$
these components either lie in the space $x_1\neq0$ or intersect the hyperplane $x_1=0$;

{\rm a)} the components lying in the space $x_1\neq0$ form pairs consisting of manifolds symmetric with respect to the reflection $x_1\mapsto-x_1$; these pairs are numbered by the locations of zero and all real roots of the polynomial {\rm(\ref{P1})} that has $k_i$ multiple real roots of multiplicity $\mu_i+1,i=1,\dots,p$ and does not satisfy the condition {\rm(\ref{cond2})};

{\rm b)} the components intersecting the hyperplane $x_1=0$ are numbered by the locations of all real roots of the polynomial {\rm(\ref{P1})} that has $k_i$ nonzero multiple real roots of multiplicity $\mu_i+1,i=1,\dots,p$ and satisfies the condition {\rm(\ref{cond2})};

{\rm2)} if $\mathcal{A}$ is divisible by $A_1^2$, then there are components lying in the space
$$
\Pi_2: x_1=u=0,x_2<0;
$$
these components are numbered by the locations of all real roots of the polynomial {\rm(\ref{P1})} that has $m_i$ nonzero multiple real roots of multiplicity $\nu_i+1,i=1,\dots,r$, where $A_{\nu_1}^{m_1}\dots A_{\nu_r}^{m_r}A_1^2=\mathcal{A}$;

{\rm3)} if $\mathcal{A}$ is divisible by $A_3$, then there are components lying in the space
$$
\Pi_3: x_1=x_2=u=0;
$$
these components are numbered by the locations of all real roots of the polynomial {\rm(\ref{P1})} that has $m_i$ nonzero multiple real roots of multiplicity $\nu_i+1,i=1,\dots,r$, where $A_{\nu_1}^{m_1}\dots A_{\nu_r}^{m_r}A_3=\mathcal{A}$.}

\medskip
Theorem 1 reduces the problem of calculating the adjacency indices $J_{A_{\mu_1}^{k_1}\dots A_{\mu_p}^{k_p}}(D_\mu^\delta)$ to the combinatorics of possible locations of real roots of the polynomial (\ref{P1}) on $t_2$-axis taking into account their multiplicities and the position of zero. This combinatorics is described below. 

Let
$$
\left\langle n_1,\dots,n_r\right\rangle=\frac{(n_1+\dots+n_r)!}{n_1!\dots n_r!}
$$
for any nonnegative integers $n_1,\dots ,n_r$. If necessary, we will write the numbers $n_1,\dots ,n_r$ in the expression $\left\langle n_1,\dots,n_r\right\rangle$ in two lines:
$$
\left\langle n_1,\dots,n_r\right\rangle\equiv\left\langle n_1,\dots ,n_p\atop n_{p+1},\dots ,n_r\right\rangle.
$$

\medskip
{\bf Theorem 2.} {\it Let $\mathcal{A}=A_{\alpha_1}^{i_1}\dots A_{\alpha_k}^{i_k}A_{\beta_1}^{j_1}\dots A_{\beta_l}^{j_l}$, where $\alpha_1,\dots,\alpha_k$ are pairwise distinct even numbers, and $\beta_1,\dots,\beta_l$ are pairwise distinct odd numbers. By $a_\mu$ denote the number of factors in $\mathcal{A}$ that are equal to $A_\mu$. Let
$$
m_1=\sum\nolimits_{p=1}^{k}i_p;\quad
m_2=\sum\nolimits_{r=1}^{l}j_r;
$$
$$
N=\mu-\sum\nolimits_{p=1}^ki_p(\alpha_p+1)-\sum\nolimits_{r=1}^lj_r(\beta_r+1);\quad 
\Delta(i_0)=\left[\frac{3-\delta+2(i_0+m_1)}{4}\right].
$$
Then
\begin{equation}
J_{\mathcal{A}}(D_\mu^\delta)=I_0+I_1+a_1(a_1-1)I_2+a_3I_3,
\label{main}
\end{equation}
where
$$
I_0=\sum_{0\leq i_0\leq N,\atop i_0\equiv N\,(\mathrm{mod}\,2)}\left\langle i_0,\dots ,i_k\atop j_1,\dots ,j_l\right\rangle \left(2\Delta(i_0)\frac{1+i_0+m_1+m_2}{1+i_0+m_1}-i_0\right);
$$
$$
I_1=\sum_{0\leq i_0\leq N-2\atop i_0\equiv N\,(\mathrm{mod}\,2)} \left\langle i_0,\dots ,i_k\atop j_1,\dots ,j_l\right\rangle \frac{1+i_0+m_1+m_2}{1+i_0+m_1}\left(1+i_0+m_1-\Delta(i_0)\right);
$$
\medskip
$$
I_2=\begin{cases}
0&\text{if $a_1\leq1$},\\
\sum\limits_{0\leq i_0\leq N+2\atop i_0\equiv N\,(\mathrm{mod}2)}\left\langle i_0,\dots ,i_k\atop j_1,\dots ,j_l\right\rangle
\frac{\Delta(i_0)}{(1+i_0+m_1)\,(i_0+m_1+m_2)}&\text{if $a_1\geq2$};
\end{cases}
$$
$$
I_3=\begin{cases}
0&\text{if $a_3=0$},\\
\sum\limits_{0\leq i_0\leq N+1\atop i_0\equiv N+1\,(\mathrm{mod}2)}
\left\langle i_0,\dots ,i_k\atop j_1,\dots ,j_l\right\rangle
&\text{if $a_3\geq1$}.
\end{cases}
$$
}

\medskip
{\sc Proof.} Let the map $f$ have a multisingularity of type $\mathcal{A}$ at a point $y=(x,q,u)$. Consider the polynomial (\ref{P1}) corresponding to $y$. Its real roots of odd multiplicity divide $t_2$-axis into (open) intervals. An interval is called negative (positive) if the polynomial is negative (positive, respectively) everywhere in this interval, except for roots of even multiplicity.

Let $i_0$ be the number of simple real roots of the polynomial (\ref{P1}). Then: 

$i_0\equiv N\,(\mathrm{mod}\,2)$ and $0\leq i_0\leq N$ if $y\in\Pi_1$; 

$i_0\equiv N\,(\mathrm{mod}\,2)$ and $0\leq i_0\leq N+2$ if $y\in\Pi_2$; 

$i_0\equiv N+1\,(\mathrm{mod}\,2)$ and $0\leq i_0\leq N+1$ if $y\in\Pi_3$.

Assume that $x_1\neq0$. Then the polynomial (\ref{P1}) has $m_1$ geometrically different multiple real roots of odd multiplicity and $m_2$ geometrically different roots of even multiplicity. Namely, it has $i_p$ roots of odd multiplicity $\alpha_p+1$ for each $p=1,\dots,k$ and $j_r$ roots of even multiplicity $\beta_r+1$ for each $r=1,\dots,l$. All of them are nonzero.

Simple and multiple real roots of odd multiplicity can be placed on $t_2$-axis in
\begin{equation}
\left\langle i_0,\dots,i_k\right\rangle
\label{odd-k}
\end{equation}
ways. These roots divide the axis into $1+i_0+m_1$ intervals. Zero belongs to a negative interval. The numbers of negative and positive intervals are equal to $\Delta(i_0)$ and $1+i_0+m_1-\Delta(i_0)$, respectively. The number of all locations of real roots of even multiplicity of the polynomial (\ref{P1}) on $t_2$-axis at a fixed position of zero and roots of odd multiplicity is equal to
$$
\left\langle 1+i_0+m_1,j_1,\dots,j_l\right\rangle.
$$

Thus by Theorem 1, it follows that the number of connected components of the manifold $\mathcal{A}_f$ lying in $\Pi_1$ and intersecting the space $x_1>0$ is equal to
$$
C_0=\sum_{0\leq i_0\leq N,\atop i_0\equiv N\,(\mathrm{mod}\,2)}\left\langle i_0,\dots,i_k\right\rangle \left\langle 1+i_0+m_1,j_1,\dots,j_l\right\rangle \Delta(i_0).
$$
The number of connected components of the manifold $\mathcal{A}_f$ lying in $\Pi_1$ and intersecting the space $x_1=0,u\neq0$ is equal to
$$
C_1=\sum_{0\leq i_0\leq N,\atop i_0\equiv N\,(\mathrm{mod}\,2)}\left\langle i_0,\dots,i_k,j_1,\dots,j_l\right\rangle i_0.
$$
The number of connected components of the manifold $\mathcal{A}_f$ lying in $\Pi_1$ and intersecting the space $x_1=0,u=0,x_2>0$ is equal to
$$
I_1=\sum_{0\leq i_0\leq N-2,\atop i_0\equiv N\,(\mathrm{mod}\,2)}\left\langle i_0,\dots,i_k\right\rangle \left\langle 1+i_0+m_1,j_1,\dots,j_l\right\rangle \left(1+i_0+m_1-\Delta(i_0)\right).
$$
According to the inclusion-exclusion principle, the number of all connected components of the manifold $\mathcal{A}_f$ lying in the space $\Pi_1$ is equal to $I_0+I_1$, where $I_0=2C_0-C_1$.

It remains to calculate the number of connected components of the manifold
$\mathcal{A}_f$ lying in the spaces $\Pi_2$ and $\Pi_3$. The number of locations of all real roots of the polynomial (\ref{P1}) on $t_2$-axis for $x_1=u=0,x_2<0$ is equal to the product of the numbers (\ref{odd-k}), $\Delta(i_0)$ and
$$
\frac{(i_0+m_1+m_2-1)!}{(1+i_0+m_1)!\,j_1!\dots j_l!}a_1(a_1-1).
$$
The last one is the number of locations of real roots of even multiplicity. Thus the number of connected components of the manifold $\mathcal{A}_f$ lying in $\Pi_2$ is equal to $a_1(a_1-1)I_2$.

The number of locations of all real roots of the polynomial (\ref{P1}) for $x_1=x_2=u=0$ is equal to
$$
\frac{(i_0+m_1+m_2)!}{i_0!\,i_1!\dots i_k!\,j_1!\dots j_l!}a_3.
$$
Therefore the number of connected components of the manifold $\mathcal{A}_f$ lying in $\Pi_3$ is equal to $a_3I_3$. This completes the proof of Theorem 2.

\medskip
{\bf Remark 1.} The calculation of the adjacency indices of Legendrian monosingularities of types $D_\mu^\pm,\mu\leq6$ to multisingularities of types $A_{\mu_1}^{k_1}\dots A_{\mu_p}^{k_p}$ using the formula (\ref{main}) gives the same results as \cite[Theorem 2.8]{Sed2012}.

\medskip
{\bf Remark 2.} Using (\ref{main}), it is easy to obtain the formulas
$$
J_{\bf 1}(D_{2k}^\delta)=\frac{k(k+1)}{2}+(k+2)\left[\frac{3-\delta}{4}\right],\quad J_{\bf 1}(D_{2k+1})=\frac{(k+1)(k+2)}{2},
$$
where $k\geq2$. These formulas were obtained previously by Vassiliev in \cite{Vas-discr} when studying the topology of the subsets $S(t,q)+x_1t_1+x_2t_2\leq u$ in the plane $t=(t_1,t_2)$ depending on the parameters $q,x_1,x_2,u$.

\medskip 
The author thanks M.~E.~Kazarian for helpful discussions.

\medskip

\begin{flushleft}

Gubkin University, Moscow, Russia

\medskip

vdsedykh@gmail.com

sedykh@mccme.ru

\end{flushleft}


\newpage

\begin{thebibliography}{9}

\bibitem{Arn96}
V. I. Arnold, {\it Singularities of caustics and wave fronts}, Math. Appl. (Soviet Ser.), vol. 62, Kluwer Acad. Publishers, Dordrecht 1990.

\bibitem{Loenga}
E. Looijenga, {\it The discriminant of a real simple singularity}. Compositio Math. {\bf 37} (1978), no. 1, 51--62.

\bibitem{Sed2012}
V. D. Sedykh, {\it On the topology of wave fronts in spaces of low dimension}. Izv. Math.
{\bf 76} (2012), no.~2, 375--418.

\bibitem{SedARMJ}
V. D. Sedykh, {\it On Lagrangian and Legendrian singularities}. Arnold Math. J. {\bf 7} (2021), no.~2, 195--212.

\bibitem{Sed-kniga-2021}
V. D. Sedykh, {\it Mathematical methods of catastrophe theory}. Moscow Center of Continuous Mathematical Education, Moscow, 2021, 224 pp. (Russian).

\bibitem{Sed2024}
V. D. Sedykh, {\it On the adjacency of type $D$ singularities of a front.} Russian Math. Surveys, {\bf 79} (2024), no.~3, 550--552.

\bibitem{Vas-discr}
V. A. Vassiliev, {\it Complements of discriminants of simple real function singularities}. Israel J. of Math. (2024). https://doi.org/10.1007/s11856-024-2627-8.

\end{thebibliography}
\end{document}